\documentclass[11pt,leqno]{article}
\usepackage{graphicx, amsfonts, amsthm, amsxtra, verbatim, multicol}
\usepackage[mathscr]{euscript}
\newtheorem{theo}{Theorem}

\newtheorem{coro}{Corollary}

\newtheorem{lemm}{Lemma}
\newtheorem{prop}{Proposition}
\newtheorem{rem}{Remark}

\def\End{{\mathrm End}}              

\newcommand\Pn[1]{\mathbb{P}^{#1}}   
\def\Z{\mathbb{Z}}                 
\def\ZZ{\mathbb{Z}}               
\def\Q{\mathbb{Q}}                   
\def\C{\mathbb{C}}                   
\newcommand\ep[1]{e^{\frac{2\pi i}{#1}}}
\newcommand{\mat}[4]{
     \begin{pmatrix}
            #1 & #2 \\
            #3 & #4
       \end{pmatrix}
    }                                
\newcommand{\matt}[2]{
     \begin{pmatrix}                 
            #1   \\
            #2
       \end{pmatrix}
    }
\def\ker{{\mathrm ker}}              
\def\dR{{\mathrm dR}}                

\def\hc{{\mathsf H}}                 
\def\pedo{{\cal P}}                  

\def\gr{{\mathrm Gr}}                
\newcommand\SL[2]{{\mathrm SL}(#1, #2)}    

\def\ZZ{{\mathbb Z}}

\def\CC{{\mathbb C}}
\def\n{{\mathfrak{n}}}
\def\p{{\mathfrak{p}}}
\newcommand{\si} {\sigma}
\def\ord{{\rm ord}}              
\begin{document}
\begin{center}
{\LARGE\bf Hypergeometric series and Hodge cycles 
of four dimensional cubic hypersurfaces
}
\footnote{ Math. classification: 14C30, 14Kxx
} \\

\vspace{.25in} {\large {\sc H. Movasati}} \\
IMPA - Instituto de Matematica Pura e Aplicada \\
Estrada Dona Castorina 110 \\
Jardim Botanico  22460-320 \\
Rio de Janeiro  - RJ  Brasil \\
Email: {\tt hossein@impa.br} \\
{\large {\sc  S. Reiter}} \\
Technische Universit\"at Darmstadt,
Fachbereich Mathematik \\
Schlo\ss gartenstr. 7,
64289 Darmstadt, Germany \\
Email: {\tt reiter@mathematik.tu-darmstadt.de}
\begin{abstract}
In this article we find connections between
the values of Gauss hypergeometric functions and
the dimension of the vector space of Hodge cycles of four dimensional
cubic  hypersurfaces. Since the Hodge conjecture is well-known for those
varieties we calculate values of Hypergeometric
series on certain $CM$ points.
Our methods is based on the calculation
of the Picard-Fuchs equations in higher dimensions, reducing them to the
Gauss equation and then applying the Abelian Subvariety Theorem to
the corresponding  hypergeometric abelian varieties.
\end{abstract}
\end{center}
{\tiny \tableofcontents }
\section{Introduction}
Recently, there have been attempts to relate algebraic values of 
hypergeometric functions in algebraic points to certain phenomenons
in Arithmetic Algebraic Geometry. In this direction one can mention
F. Beukers article \cite{beu95} in which he calculates 
explicitly the values of (elliptic) hypergeometric functions in
some $CM$ points. Such a hypergeometric function can be rewritten as 
an elliptic integral and so the word $CM$ comes from the corresponding 
$CM$ elliptic curves. H. Shiga and J. Wolfart 
in \cite{wo04} (see also \cite{shwo95} and \cite{shtswo}) have studied
the algebraic values of Schwarz functions (certain quotient of two 
hypergeometric functions) in algebraic points. 
Their geometric tool is Abelian Subvariety Theorem and
its consequence on the dimension of the periods of an abelian variety. 
Again, the idea behind algebraicity of a  value of the Schwarz function is
that the corresponding abelian variety has more automorphism than the
neighboring varieties. For a general survey about algebraic and transcendental
periods see Waldschmidt's expository article \cite{wal}. 

The idea of this article took place in our mind when we observed that
many abelian integrals of affine hypersurfaces satisfy second order 
Picard-Fuchs differential equations and they can be written in terms of 
hypergeometric functions. A consequence of the Hodge conjecture 
(see Proposition \ref{deligne}) gives special values of such integrals 
on the so called Hodge cycles. This motivated us to use the mentioned 
geometric phenomenon and obtain special values of hypergeometric functions. 
Fortunately, the Hodge conjecture is known in all cases which we use. 
We concentrate mainly on cubic hypersurfaces and we do not yet know whether 
one can  classify the hypergeometric functions appearing in this way. 
Below is the summary of the results of this article.

Let $n\geq 2 $ be an even number.
We consider a family of $n$-dimensional varieties $M:=\{M_t\}_{t\in\Pn 1}$ and
a meromorphic differential $n$-form $\omega$ whose restriction to each fiber has
no residues  around poles. We  assume that:
\begin{enumerate}
\item
For all $t\in \Pn 1$ except a finite number  of them, the fiber $M_t$ is
smooth.
\item
The Hodge decomposition of $H^n(M_t,\C)$ is of the form
$$
H^n(M_t,\C)=H^{\frac{n}{2}-1, \frac{n}{2}+1}
\oplus H^{\frac{n}{2}, \frac{n}{2}}\oplus H^{\frac{n}{2}+1, \frac{n}{2}-1}, 
$$
$$
\dim_\C(H^{\frac{n}{2}-1, \frac{n}{2}+1})=
\dim_\C(H^{\frac{n}{2}+1, \frac{n}{2}-1})=1
$$
and $\omega\mid_ {M_t}$ form a basis of $H^{\frac{n}{2}-1, \frac{n}{2}+1}$.
\item
The Picard-Fuchs equation of $\int_{\delta_t} \omega$, where
$\delta_t\in H_n(M_t,\Z)$ is a continuous family of cycles,
 is a pull-back of
the Gauss equation
\begin{equation}
\label{gauss}
z(1-z)y''+(c-(a+b+1)z)y'-aby=0.
\end{equation}
\end{enumerate}

For a number field $k$,
the vector space of Hodge cycles with coefficients in $k$
is defined as follows:
\begin{equation}
\label{30apr05}
\hc_t(k):=H^{\frac{n}{2},\frac{n}{2}}\cap H^n(M_t,k).
\end{equation}
The usual definition of a Hodge cycles is with $k=\Q$.
Since the hypergeometric function
$$
F(a,b,c|z)=\sum_{n=0}^{\infty} \frac{(a)_n(b)_n}{(c)_n n!}z^n, \ c\not\in
\{0,-1,-2,-3,\ldots \},
$$
where  $(a)_n:=a(a+1)(a+2)\cdots (a+n-1)$,
solves the Gauss equation we find certain relations between the values
of $F$-functions and Hodge cycles of the fibers $M_t$. For instance,
the family
\begin{equation}
\label{12july05}
M_t:\ f(x):=x_1^3+x_2^3+\cdots+x_5^3-x_1-x_2=t,
\end{equation}
satisfies the conditions 1,2 and 3 and  has three
critical fibers corresponding to
$t=0,\pm 2 (\frac{2}{3})^\frac{3}{2}$. Here $\omega=\nabla \frac{dx}{df}$,
where $dx$ is the wedge product of all $dx_i$'s , $\frac{dx}{df}$ is
the Gelfand-Leary form of $dx$  and
$\nabla$ is the Gauss-Manin connection with respect to the parameter $t$
(see \S \ref{gausssection}). Let $B(a,b)=
\frac{\Gamma(a)\Gamma(b)}{\Gamma(a+b)}$ be the beta function. For two
numbers $r,s\in \C$ we say that $r\sim s$ if $\frac{r}{s}\in\bar \Q$. 
We prove that:
\begin{theo}
\label{main1}
Let $k$ be a number field with $\Q(\zeta_3)\subset k$,
\begin{equation}
\label{26.5}
E_t:\ y^2=x^3-3x+2-\frac{27}{4}t^2,\ t\in \bar\Q
\end{equation}
and $z=\frac{27}{16}t^2$.
\begin{enumerate}
\item
If $E_t$ is not a $CM$-elliptic curve then the
$k$-vector space $\hc_t(k)$ has codimension $2$ ;
\item
For those values of $t\in\bar\Q$ such that the elliptic curve $E_t$
is $CM$, the  following value of the Schwarz function
\begin{equation}
\label{26.5.05}
D(0,0,1|z):= -e^{-\pi i \frac{5}{6}} 
\frac{ F(\frac{5}{6},\frac{1}{6},1|z)}
 {  F(\frac{5}{6},\frac{1}{6},1|1-z)}
\end{equation}
is algebraic and if it belongs to $k$
then the  $k$-vector space
$\hc_t(k)$ has codimension one.
\item
The number (\ref{26.5.05}) belongs to $\Q(\zeta_3)$ for  some 
$t\in\bar\Q$ if and only if
$$
 F(\frac{5}{6},\frac{1}{6},1|z) \sim
\frac{1}{\pi^2}\Gamma(\frac{1}{3})^3,\ 
 F(\frac{5}{6},\frac{1}{6},1|1-z)\sim
\frac{1}{\pi^2}\Gamma(\frac{1}{3})^3.
$$
\end{enumerate}
\end{theo}
We have $j(E_t)=\frac{-2^43^3}{z(z-1)}$ and, for instance, if
$j(E_t)=2^43^35^3,\ -2^{15}3\cdot5^3$ then the condition of the third part of 
the above theorem is satisfied (see \cite{sil} p.483).
The proof of the first and second part uses a simple consequence of 
Abelian Subvariety Theorem on the periods of an elliptic curve. 
The proof of the third
part is based on the Hodge conjecture for the variety $M_t$. We note that 
this conjecture is proved for cubic hypersurfaces of dimension $4$ by 
C. Clemens, J. P. Murre and S. Zucker (see \cite{zu77} and its references).

Let us now explain the content of each section. In section
\ref{hodgecycles} we introduce the notion of a Hodge cycle by means of 
vanishing of abelian integrals. \S \ref{pivalue} is devoted to a consequence
of the Hodge conjecture on the values of integrals over Hodge cycles.
In \S \ref{gm} we recall the notion of Gauss-Manin 
connection associated to a fibration and state a conjecture which has been
verified for some examples.
\S \ref{chishod} is devoted to the examples of Hodge structures 
whose Hodge numbers is of the type $1,x,1$.  
In \S \ref{AST} we use Abelian Subvariety theorem and derive the fact that
a simple abelian variety $A$ is $CM$ if and only 
if the periods of a differential form of the first kind on
$A$ generate a $\bar\Q$-vector space of dimension one.
In \S \ref{gausssection} we fix up the notations related to
the Gauss-equation and hypergeometric functions. Following
the ideas of \cite{shtswo}, \cite{wo04} in \S \ref{familyofcurves} 
we introduce a one dimensional integral representation of hypergeometric 
functions.
\S \ref{four} is devoted to  the calculations related to the example
(\ref{26.5}). In this section we prove Theorem \ref{main1}. 
In the last section we give many other examples of
cubic hypersurfaces which can be analyzed by the methods of this article
and hence may result to some algebraic relations between the values
of hyper geometric functions.
\section{Families of varieties and Hodge cycles}
\label{first}
The classical definition of a Hodge cycle over $k=\Q$ is given in
(\ref{30apr05}). This definition can be formulated by means of 
$n$-dimensional
abelian integrals. Since such abelian integrals satisfy
Picard-Fuchs equations, our problem of studying Hodge cycles reduces to the
study of Picard-Fuchs equation.
For $d$ an integer we set
$$
G_d:=\{\zeta_d^i\mid i=0,1,\ldots,d-1\},\ \zeta_d=\ep d
$$
\subsection{Hodge cycles}
\label{hodgecycles}
In the following sections we will take a family $\{M_t\}_{t\in\Pn 1}$ of
$n$ dimensional varieties and compute certain Picard-Fuchs equations. In all
examples we take $f$ a polynomial of degree $d$
in $\C^{n+1}$ and set $L_t:=\{f=t\}, M_t:=\overline{L_t}$,
where the closure is taken in $\Pn {n+1}$. We assume that $M_t$ is smooth for
all $t$ except a finite number of them and denote by $k$ a number field.

A cycle $\delta\in H_n(M_t,k)$ is called Hodge if
$$
\int_{\delta}\omega=0,\ \omega\in F^{\frac{n}{2}+1}H_\dR(M_t),
$$
where
$F^\bullet$ denotes the Hodge filtration of $H_\dR^n(M_t)$. By Poincar\'e
duality this definition of a Hodge cycle translate into the definition
(\ref{30apr05}) in the Introduction.
If the support of $\delta$ is in $L_t$ then we can reformulate this
definition using the mixed Hodge structure of $L_t$, i.e. $\delta\in H_n(L_t,k)$ is Hodge
if
$$
\int_\delta\omega=0, \ \omega\in \gr^{i}_F\gr^W_n H^n_\dR(L_t), \ i\geq
\frac{n}{2}+1,
$$
where
\begin{equation}
 0=F^{n+1}\subset F^n\subset \cdots \subset F^1
\subset F^0=H^n_\dR(L_t),
\end{equation}
$$
0=W_{n-1}\subset W_n\subset W_{n+1}=H^n_\dR(L_t)
$$
are the Hodge and weight filtrations of the Mixed Hodge structure of
$H_\dR^n(L_t)$.
We denote by $\hc_t(k)$ the $k$-vector space
of Hodge cycles  in $H_n(M_t,k)$.
For more details see (\cite{mo}).

\begin{rem}\rm
For two number fields $k_1\subset k_2$ we do not have necessarily
$\hc_t(k_1)\otimes_{k_1} k_2=\hc_t(k_2)$. This can be seen in the next
section.
\end{rem}
\subsection{Hodge cycles and algebraic values of abelian integrals}
\label{pivalue}
We follow the notations of the previous section.
Let $\delta\in H_n(L_t,\Q)$ be a Hodge cycle.
The Hodge conjecture claims that  there is an algebraic cycle
$Z=\sum_{i=1}^kr_i Z_i,\ \dim_\C Z_i=\frac{n}{2}$ such that the topological
class $[Z]:=\sum_{i=1}^kr_i [Z_i]\in H_n(M_t,\Q)$ of $Z$ is equal to
the image of $\delta$ in $H_n(M_t,\Q)$.
In other words, every Hodge cycle is an algebraic cycle.
\begin{prop}
\label{deligne}
If $[Z]\in H_n(L_t,\Q)$ is an algebraic cycle and $L_t$ is defined over
$\bar \Q$ then for any differential form $\omega$ on $L_t$ which does not have residues at
infinity and is 
defined over $\bar\Q$ we have
$$
\int_{[Z]}\omega \in (2\pi i)^{\frac{n}{2}}\bar\Q.
$$
\end{prop}
See Deligne's lecture \cite{dmos} Proposition 1.5 for a proof. Note that
if $M_t$ is defined over $\bar\Q$ and $[Z]$ is an algebraic cycle
we can assume that $Z$ is defined over $\bar\Q$.
Deligne has applied the above proposition to Fermat hypersurfaces
and has obtained  algebraic relations between the values of $\Gamma$
function on rational points (see Theorem 7.15 \cite{dmos}). Our
approach to the algebraic relations between the values of $F$-functions
follows the same method used by Deligne.
\subsection{Gauss-Manin connection and Picard-Fuchs equations}
\label{gm}
In this section we consider a family of algebraic varieties
${\cal M}$ over $\C(t)$ and a differential 1-form
$\omega\in H^n_\dR({\cal M})$.  We denote by $M_t$ the specialization
of ${\cal M}$ to a fixed
$t\in \C$ and we call them fibers.
We use also the notation $\omega$ for  its specialization
to the fiber $M_t$; being clear in the text what we mean.
Let $S\subset \C$ be the locus of singular fibers.
The multi valued functions
$$
I_{\gamma_t}(t):=\int_{\gamma_t} \omega
$$
spans the solution space of  a Picard-Fuchs equation
\begin{equation}
\label{pf}
p_m(t)y^{(m)}+p_{m-1}(t)y^{(m-1)}+\cdots+p_1(t)y'+p_0(t)y=0,\ m\leq 2g,\
p_i\in \C[t],
\end{equation}
where $2g= \dim H^n_\dR(\cal M)$ and $\gamma_t\in H_n(M_t,\Z)$ is a
continuous family of cycles which can be  defined in any simply connected
region in $\C\backslash S$. Moreover, if
the degree $m$ of (\ref{pf}) is equal to $2g$ then a basis
$\gamma_{i,t}\in H_n(M_t,\Z),\ i=1,2,\ldots, 2g$ gives a basis
$I_{\gamma_{i,t}}, \ i=1,2,\ldots,2g$
of (\ref{pf}). One obtains (\ref{pf}) in the
following way:
 The Gauss-Manin connection $\nabla$ is defined in the cohomology
$H^n_\dR({\cal M})$ and it has poles in $S$. Its main property is
\begin{equation}
\label{20.1}
\frac{d}{dt}\int_{\gamma_t} \omega= \int_{\gamma_t} \nabla \omega.
\end{equation}
Since $H^n_\dR({\cal M})$ is a $\C(t)$-vector space of dimension $2g$,
there must be  $\C[t]$-linear relation between
$$
\omega, \nabla\omega,\ldots, \nabla^{2g} \omega,
$$
integrating such linear equation and using (\ref{20.1}) we get
(\ref{pf}). 
The set $S':=\{t\in \C \mid p_m(t)=0\}$ is called the singular set of
(\ref{pf}) and we know that $S\subset S'$.
The points in $S'\backslash S$ are called apparent singularities
(see \cite{ano}). Let 
\begin{equation}
\label{22june2005}
\pedo({\cal M},\omega):=\{ \int_{\gamma_{t}} \omega 
\mid \gamma_{t}\in H_{n}(M_{t},\Q)\}.
\end{equation}
In the next sections we will need to verify the following conjecture
for our example.

{\bf Conjecture:}
Let $({\cal M}_i,\omega_i),\ i=1,2$ be two pairs as above with 
$\dim {\cal M}_1\geq \dim {\cal M}_2$. If the associated
period space $\pedo({\cal M}_i,\omega_i),\ i=1,2$
span the solution space of the same Picard-Fuchs equation
then there is a non-zero complex number $\gamma$ and a number field $k$
such that 
$\pedo({\cal M}_1,\omega_1)$ is generated over $\Q$ by $\gamma\cdot k \cdot 
\pedo ({\cal M}_2,\omega_2)$.

One may expect that the number $\gamma$ and the number field $k$ depend 
only on some numerical invariants of the pairs 
$({\cal M}_i,\omega_i),\ i=1,2$. 

In the examples which we have verified this conjecture, the
number $\gamma$ is obtained by the values of the beta function on rational 
points. 
\subsection{Four dimensional cubic hypersurfaces}
\label{chishod}
In this section we consider the polynomial
$f=x_1^d+x_2^d+\cdots+x_{n+1}^d-g(x),\ \deg(g)<d $ and the associated
family of varieties
$$
M_t:\overline{\{x\in\C^{n+1} \mid f(x)=t\}}.
$$
In this article we only deal with the following cases:
$$
 d=3,\ n=4,\ h^{0,4}=h^{4,0}=0,\ h^{13}=h^{31}=1, \ h^{2,2}=20,
$$
where $h^{i,j}$'s are the Hodge numbers of a smooth $M_t$ in the corresponding
case (see for instance \cite{kk98}).
The $5$-form $\omega=Pdx$, where
$dx:=dx_1\wedge dx_2\wedge\cdots\wedge dx_{n+1}$
and $P\in\C[x]$,  can be interpreted as a section of the cohomology bundle of
the fibration $\{M_t\}_{t\in \Pn 1}$ by considering the Gelfand-Leray form
$\frac{\omega}{df}$ (see for instance \cite{arn}).
A way of constructing such section is as follows:
We take the residue of $\frac{\omega}{f-t}$ on $M_t$, which is an element
in $H^n(M_t,\C)$.
\begin{theo}
\label{mova}
If $n=4,d=3$  then
the differential form $\nabla(\frac{dx}{df})$ is a basis of $H^{3,1}$, where
$H^4_\dR(M_t)=H^{3,1}\oplus H^{2,2}\oplus H^{1,3}$
 is the Hodge
decomposition.
\end{theo}
The above theorem for the residue of $\frac{dx}{f-t}$ in $M_t$ is a particular
case of Griffiths theorem for the Hodge filtration of hypersurfaces
(see \cite{gr69}, \cite{kk98} and also \cite{mo} Theorem 2).
The fact that the mentioned residue can be expressed by the Gauss-Manin
connection of $\frac{dx}{df}$ is explained in Lemma 2 \cite{mo}. See also
this article for a generalization of the above theorem. 
\subsection{Abelian Subvariety Theorem}
\label{AST}
In this section we use Abelian Subvariety Theorem and
calculate the dimension of the $\bar\Q$-vector spaces spanned by
the periods of a differential 1-form of the first kind on a 
$CM$ abelian variety.

Let $A$ be an abelian variety defined over $\bar \Q$ and 
$t_A$ be the tangent space at the origin $0\in A$. The exponential map
$$
\exp :t_A(\C)\rightarrow A(\C)
$$
is well-defined and set $\Lambda:=\exp^{-1}(0)$. Here, by $t_A(\C)$ we mean the
complexification of $t_A$.
Let $W$ be a proper linear subspace of $t_A$
(defined over $\bar \Q$) such that
\begin{equation}
\label{peri}
0\neq \gamma\in \Lambda \cap W(\C).
\end{equation}
\begin{theo}
Under the above conditions, there exists a proper connected abelian sub
variety $B\subset A$ (defined over $\bar \Q$) such that
$$
t_B\subset W,\ \hbox{ and } \gamma\in t_B(\C).
$$
\end{theo}
For this version of Abelian Subvariety Theorem see for instance
Lemma 1 of \cite{shwo95}. 

\begin{coro}
\label{bewerbung}
For a simple abelian variety  $A$ and a differential form of the first kind
$\omega$, both defined over $\bar \Q$, the $\bar\Q$-vector space 
 $\{\int_\delta\omega \mid  \delta\in H_1(A(\C),\bar \Q)\}$ is of 
dimension one if and only if  $A$ is $CM$ and $\omega$ is an eigendifferential under
the action of the $CM$ field.
\end{coro}
\begin{proof}
This lemma is  Proposition 3 of \cite{shwo95}.
Let $g$ be the dimension of $A$. 
We take a $\Q$-basis $\delta_i, i=1,2,\ldots,2g$ of
the homology $H_1(A(\C),\Q)$. Without losing generality, we 
assume that $\int_{\delta_1}\omega\not =0$. 
The non trivial part of the theorem is to prove that if there are 
constants $a_i\in \bar \Q,\ i=2,3,\ldots,2g$ such that 
\begin{equation}
\label{1jun05}
\int_{\delta_i}\omega+a_i\int_{\delta_1}\omega=0,\ i=2,3,\ldots,2g 
\end{equation}
then $A$ is $CM$. To prove this we define 
$$
B=A\times A,\ \omega_i=\omega \oplus (a_i\omega)\in \Omega_{B}, \ 
\delta=\delta_i\oplus \delta_1\in H_1(B(\C),\Q)
$$
for a fixed $i$, 
where the the sums are well-defined 
using the isomorphisms $\Omega_B\cong\Omega_A\oplus\Omega_A$ and 
$H_1(B,\Q)\cong
H_1(A,\Q)\oplus H_1(A,\Q)$. 
Now the condition (\ref{1jun05}) implies that
$$
\int_\delta \omega'=0.
$$
Looking $\omega'_0$ as a linear map from $t_A$ to $\bar \Q$ this implies that
$\delta\in W:=\ker(\omega'_0)$. Therefore, by Abelian Subvariety Theorem there
exists an abelian subvariety $B'\subset B$ such that $\delta$ is 
supported in  $H_1(B',\Q)$ and the restriction of $\omega'$ to $B'$ is 
identically zero. Since $A$ is simple, $B'$ is of dimension $g$ and
the projections $\pi_i, i=1,2$ are isogeny between $B'$ and $A$ and 
$\pi_2\circ \pi_1^{-1}$ gives a non trivial endomorphism of $A$ because it
sends $\delta_1$ to $n\delta_i$ for some $n\in\Z$.
We have proved that $H_1(A(\C),\Q)$ as $\End_0(A)$-module is of dimension
one which implies that $A$ is $CM$.   
\end{proof}
\section{Hypergeometric functions}
Recall the definition of Pochhammer cycles and  the corresponding
Kummer solutions form \cite{IKSY} and \cite{shtswo} \S 5. The Pochhammer
cycle associated to the points $a,b\in\C$ and a path $s:[0,1]\rightarrow \C$ 
connecting $a$ to $b$, is the commutator
$$
[\gamma_b,\gamma_a]=\gamma_b^{-1}\cdot \gamma_a^{-1}\cdot \gamma_b\cdot 
\gamma_a,
$$
where $\gamma_b$ is a loop along $s$ starting and ending at the point 
$b=s(\frac{1}{2})$ which encircles $b$ once anticlockwise, 
and $\gamma_a$ is a similar  loop with respect to $a$. 
We will need the following simple proposition:
\begin{prop}
\label{2june}
Let $f$ be a holomorphic multi-valued function in $\C$ and $\mu,\alpha\not 
\in\Z$  be the exponents of $f$ at $z=a$ (resp. $z=b$), i.e. in a neighborhood 
of $a$ we can write
$f=g\cdot (z-a)^\mu$, $g$ a holomorphic function around $a$. 
Then for a path $s$ connecting $a$ to $b$ outside the branching points of 
$f$ we have:
$$
\int_{[\gamma_b,\gamma_a]}
f(x)dx=
( 1-e^{2\pi i \alpha})
( 1-e^{2\pi i \mu})
\int_s f(x)dx.
$$ 
\end{prop}
\subsection{Gaussian System}
\label{gausssection}
To the Gauss-equation (\ref{gauss}) we can associate the system
\begin{equation}
\label{gsystem}
 Y' =\left ( \frac{1}{z} \left(\begin{array}{cc}
                                c-1 & -b \\
                                    0 & 0
                                 \end{array}\right)
 + \frac{1}{z-1} \left(\begin{array}{cc}
                                0  & 0 \\
                                  a & c-a-b-1
                                 \end{array}\right) \right )Y, \;
\end{equation}
with the fundamental system
\begin{eqnarray*}
Y &= &
\mat
{\int_{[\gamma_0,\gamma_z]} \varphi(x,z) \frac{dx}{x} }
{\int_{[\gamma_1,\gamma_z]} \varphi(x,z) \frac{dx}{x} }
{\int_{[\gamma_0,\gamma_z]} \varphi(x,z) \frac{dx}{1-x} }
{\int_{[\gamma_1,\gamma_z]} \varphi(x,z) \frac{dx}{1-x} } \\ &= &
\mat
{\int_0^z  \varphi(x,z) \frac{dx}{x} }
{\int_1^z \varphi(x,z)  \frac{dx}{x} }
{\int_0^z \varphi(x,z)  \frac{dx}{1-x} }
{\int_1^z  \varphi(x,z) \frac{dx}{1-x}  }C=
\end{eqnarray*}
\begin{equation}
\label{akh}
\mat
{
F(a,b,c|z)
}
{
(1-z)F(1-a,1-b,1+c-a-b|1-z)
}
{
\frac{a}{c} z F(a+1,b+1,c+1|z)
}
{
\frac{c-a-b}{-b}F(-a,-b,c-a-b|1-z)
}
\end{equation}
$$
\cdot
\mat
{B(a,c-a) z^{c-1} }{0}{0}{-e^{-\pi i (1+a-c)}B(1-b,c-a)(1-z)^{c-a-b-1}}C,
$$
where
$$
\varphi(x,z):=x^{a} (1-x)^{-b} (z-x)^{c-a-1}
$$
and
$$
C:=
\mat
{(1-e^{2\pi i a })(1-e^{2\pi i(c-a-1)})}{0}{0}{(1-e^{-2\pi i b })
(1-e^{2\pi i(c-a-1)})}.
$$
The first  equality is derived from Proposition \ref{2june}.
Other equalities are obtained from Theorem 4.4.3 p. 99 together with 
chapter 1.3.  of \cite{IKSY} 
(see also \cite{shtswo} Theorem 5.3
for corrected ones). 
In particular the last equality is well-defined
when $c,c-a-b\not\in\{0,-1,-2,-3,\ldots\}$. 
Our $abc$-notations is different form those in \cite{shtswo}. Let
$a',b',c'$ be those in \cite{shtswo} p. 648. Then
$$
a=a'-c'+1,\ b=b'-c'+1,\ c=2-c' \hbox{ and }\  a'=a-c+1, \ b'=b-c+1,\  c'=2-c
$$
The corresponding angular parameters are given by 
\begin{equation}
\label{angular}
\nu_0=c-1,\ \nu_1=c-a-b,\ \nu_\infty=a-b
\end{equation}
Let us remark that the above fundamental system can also be written as
 \begin{eqnarray*}
Y & = &
\mat{z^{c-1}}{0}{0}{-b^{-1}z^c}\mat {y_1}{y_2}{y_1'}{y_2'},
\end{eqnarray*}
where
 $(y_1,y_2)$ is a basis of the solution space of (\ref{gauss}). 
For a differential system $Y'=AY$ the determinant of the 
fundamental system satisfies $\det(Y)'={\rm trace}(A)\cdot \det(Y)$ and
so in our case we have
\begin{equation}
\label{13july05}
\det(Y)\sim \pi\frac{B(a,c-a)}{B(b,c-b)} z^{c-1}(z-1)^{c-a-b-1}.
\end{equation}
The constant in the above relation is obtained by the formula
$$
F(a,b,c\mid 1)=\frac{\Gamma(c)\Gamma(c-a-b)}{\Gamma(c-a)\Gamma(c-b)}
$$
and substituting $z=1$ in (\ref{akh}). 

The system (\ref{gsystem}) has monodromy
$$
A_0 =\mat
{e^{ 2\pi i c}}
{e^{-2\pi i b}-1}01 ,\
A_1:=\mat
10
{e^{2\pi i (c-a)} (e^{-2\pi i a}-1)}
{e^{2\pi i (c-a-b)}},\
A_\infty:= A_0 A_1,
$$
at $0,1, \infty$ respectively. The Schwarz function associated to
the Gauss equation (\ref{gauss}) is
\begin{equation}
\label{12ju}
D(\nu_0,\nu_1,\nu_\infty|z)=
\frac{\int_{[\gamma_0,\gamma_z]} \varphi(x,z) \frac{dx}{x}}{\int_{[\gamma_1,\gamma_z]} \varphi(x,z) \frac{dx}{x}} =
\end{equation}
$$\frac{ (1-e^{2\pi i a }) }{ {-e^{\pi i (1+a-c)} (1-e^{-2\pi i b })}} \;
  \frac{B(a,c-a)}{B(1-b,c-a)} \;
  \frac{z^{c-1}F(a,b,c|z)}{F(1-a,1-b,c-a-b+1|1-z)}.
$$
In particular, we get for $1=c=a+b$
\[ D(0, 0, 2a-1|z)= -e^{-\pi i a} \frac{ F(a,b,1|z)}
 {  F(b,a,1|1-z)} \]
and for $z=\frac{1}{2}$ we even get
 \[ D(0,0,2a-1|1/2)=-e^{-\pi i a}.\]

Sometimes it is useful to use other Pochhammer cycles instead of those
used in $Y$. The corresponding integrals are calculated in Theorem 5.3 of
\cite{shtswo} p. 648. For instance for the case $1=c=a+b$ the $(2, 2)$-entry
of $Y$ is not well-defined. We use the Pochhammer cycle
$[\gamma_0,\gamma_1]$ and for the second row of $Y$ we obtain
\begin{equation}
\label{ukh}
\matt {\int_{[\gamma_0,\gamma_1]}\varphi(x,z)\frac{dx}{x}}
{\int_{[\gamma_0,\gamma_1]}\varphi(x,z)\frac{dx}{1-x}}
\end{equation}
$$
=
(1-e^{2\pi i a})
( 1-e^{2\pi i b}) 
B(a,-b+1)z^{c-a-1}
\matt{ F(a-c+1,a,a-b+1|\frac{1}{z})}
{\frac{a}{a-b}  F(a-c+1,a+1,a-b+1|\frac{1}{z})}.
$$
Note that if $(Y_1,Y_2)^t$ is a column of $Y$ or the above matrix then
\begin{equation}
\label{khyanat}
Y_2=-b^{-1}(zY_1'+(1-c)Y_1).
\end{equation}
\subsection{A family of curves associated to the Gaussian System}
\label{familyofcurves}
Let
$$
\mu_i=a_i+b_i,\ a_i\in\Z,\ 0 \leq b_i<1
$$
and $k$ be the least common denominator of $b_0, b_1,b$. In this
section we assume that $\mu_i$'s are rational non-integer numbers. 
Define $X(k,z)$ as
the smooth curve obtained by the desingularization of
$$
y^k=x^{kb_0}(1-x)^{kb_1}(z-x)^{kb}.
$$
Define the differential forms
$$
\eta_1:=x^{-a_0}(1-x)^{-a_1}(z-x)^{-a} \frac{dx}{y},\
\eta_2:=x^{-a_0}(1-x)^{-a_1}(z-x)^{-a} \frac{xdx}{(x-1)y}.
$$
The Pochhammer cycles $[\gamma_i,\gamma_z],\ i=0,1$ lift to cycles
$\delta_i\in H_1(X(k,z),\Z)$ and
$$
Y=\mat{\int_{\delta_1}\eta_1}{\int_{\delta_2}\eta_1}
      {\int_{\delta_1}\eta_2}{\int_{\delta_2}\eta_2}.
$$
\begin{prop}
(\cite{wo88} Satz 1,2)
We have the following isogeny for the Jacobian $J_{X(k,z)}$
$$
J_{X(k,z)}\sim T(k,z)\oplus \sum_{d|k}J_{X(d,z)},
$$
where $T(k,z)$ is an abelian variety of dimension $\varphi (k)$.
Moreover
$$
\Q(\zeta_k)\subset \End_0(T):=\Q\otimes \End(T)
$$
and $H_1(T(k,t),\Q)$ as a $\Q(\zeta_k)$-vector space is of dimension two.
\end{prop}
\section{A family of cubic four dimensional Hypersurfaces}
\label{four}

The calculation of Gauss-Manin connection in higher dimensions
does not seem to be done by hand easily. In \cite{mo, hos005} the first author
has developed algorithms which calculate the Gauss Manin connection and
hence Picard-Fuchs equations  for families of hypersurfaces in weighted
projective spaces. They are implemented  in {\sc Singular} (see \cite{GPS01})
and the corresponding library is {\tt foliation.lib}.
In this and the next sections we have used this library for our calculations.
\subsection{Some Picard-Fuchs equations}
Recall the notations introduced in \S \ref{chishod} and let
$$
f=x_1^3+x_2^3+\ldots+x_5^3-x_1-x_2.
$$
The singular values of $f$ are the roots of $27t^3-16t$. Let
$$
\omega_i=\frac{x_idx}{df},\ \omega_{i,j}=\frac{x_ix_jdx}{df},\ \omega_0=\frac{dx}{df}
$$
and $\omega=(\omega_{12},\omega_2,\omega_1, \omega_0)^t$. Then
$$
\nabla(\omega)=\frac{1}{27t^3-16t}
\left ( \begin{array}{cccc}
36t^2-32/3 & -6t &  -6t & 16/9 \\
-24t & 27t^2-8 & 8 & -4t \\
-24t & 8 & 27t^2-8 & -4t \\
32 & -18t & -18t & 18t^2-16/3
\end{array} \right )\omega.
$$
This gives us a Fuchsian system of order $4$.
One gets that $\nabla\omega_0$
satisfies
\begin{equation}
\label{17Jan}
(27t^3-16t)y''+(81t^2-16)y'+15ty=0.
\end{equation}
According to Theorem \ref{mova} the differential form
$\nabla(\omega_0)$ is a basis of
$H^{31}$ and so $\delta\in H_4(M_t,\Q)$ is Hodge if and
only if
$$
\int_{\delta}\nabla\omega_0=0.
$$
\begin{rem}\rm
The form $\nabla\omega_{12}$ satisfies
$$
(27t^3-16t)y''+(81t^2-16)y'-21ty=0.
$$
Both $\omega_i, i=1,2$ satisfy
$$
(27t^3-16t)y'''+54t^2y''-3ty'+3y=0.
$$
The  $4$-forms $\omega_i,\ i=3,4,5$ satisfy
$$
9y+
(-9t)y'+
(162t^2-32)y''+
(81t^3-48t)y'''=0
$$
and $\omega_{i,j}, i=1,2, j=3,4,5$
$$
0y+
(-45t)y'+
(621t^2+32)y''+
(1134t^3-192t)y'''+
(243t^4-144t^2)y''''=0
$$
and $\omega_{ij},\ i,j=3,4,5$
$$
0y+
(-9t)y'+
(81t^2-16)y''+
(81t^3-48t)y'''=0.
$$
\end{rem}
\begin{rem}\rm
\label{29may05}
Each variety $M_t$ has the following automorphisms
\begin{equation}
\label{sarakh}
(x_1,x_2,x_3,x_4,x_5)\rightarrow (x_1,x_2,\zeta_3^{i_3}x_3, \
\zeta_3^{i_4}x_4,\zeta_3^{i_5}x_5),\ \zeta_3^{i_j}\in G_3.
\end{equation}
The family $\{M_t\}_{t\in \Pn 1}$ itself has more
automorphisms
$$
(x_1,x_2,x_3,x_4,x_5)\rightarrow (-x_1,-x_2,\zeta_6^{i_3}x_3, \
\zeta_6^{i_4}x_4, x_5\zeta_3^{i_5}),\ \zeta_6^{i_j}\in G_6,\ i_j \hbox{ odd }.
$$
Since $\nabla(\omega_0)$ is an eigen vector with eigen values in 
$\Q(\zeta_6)=\Q(\zeta_3)$ for the above automorphisms,
the period set  $\pedo(M_t,\nabla(\omega_0)),\ t\in \C\backslash C$
is a $\Q(\zeta_6)$-vector space under the usual multiplication of numbers.
\end{rem}

\subsection{Associated Gauss system}
The differential equation (\ref{17Jan})
can be written as a system
 \[ Y'= \left( \begin{array}{cc}
                      0 & 1 \\
                      -\frac{15}{27t^2-16}& - \frac{81t^2-16}{27t^3-16t}
                  \end{array}\right) Y, \]
where $Y=(y,y')^{tr}$.
Via the gauge transformation
 \[ Y \mapsto \left( \begin{array}{ll}
                      1 & 0 \\
                      0 & -3t
                  \end{array}\right)Y\]
 we obtain the Fuchsian system
 \[ Y'= \frac{1}{t} \left( \begin{array}{cc}
                      0 & -1/3 \\
                      0 & 0
                  \end{array}\right) +
           \frac{54 t }{27t^2-16}    \left( \begin{array}{cc}
                      0 & 0 \\
                      5/6 & -1
                  \end{array}\right) Y .\]
 This Fuchsian system is a rational pull back via
$z=\frac{27}{16}t^2$
 of the hypergeometric differential equation:
 \[ Y'= \frac{1}{z} \left( \begin{array}{cc}
                      0 & -1/6 \\
                      0 & 0
                  \end{array}\right) +
        \frac{1}{z-1} \left( \begin{array}{cc}
                      0 & 0 \\
                      5/6 & -1
                  \end{array}\right) Y \]
with
\begin{equation}
\label{27.5.05}
\mu_0=\mu_1=\frac{1}{6},\ \mu=\frac{5}{6},\ c=1,\ a= \frac{5}{6},\
b=\frac{1}{6}.
\end{equation}
Further,
  \[     Y   \left( \begin{array}{cc}
                      1 & 0 \\
                      0 & (1-\zeta_6)^{-1}
                  \end{array}\right) \]
has  monodromy
 \[  A_0=\left( \begin{array}{cc}
                      1 & -1 \\
                      0 & 1
                  \end{array}\right),\
 A_1=\left( \begin{array}{cc}
                      1 & 0 \\
                      1 & 1
                  \end{array}\right) \]
at $0$, respectively $1$,
and at $\infty$ it is an element of order $6$.
Thus the generated group is $\SL 2\ZZ.$

Now let us discuss the one dimensional geometric model of (\ref{17Jan}).
In this example
$$
X(6,z): y^6=x(1-x)(z-x)^{5}
$$
which is of genus $5$.
$$
\eta_1=\frac{dx}{y},\ \eta_2=\frac{xdx}{(1-x)y}
$$
$$
X(2,z): y^2=x(1-x)(z-x)^{5},\ X(3,z): y^3=x(1-x)(z-x)^{5}.
$$
With transformations $(x,y)\mapsto (x,\frac{y}{(z-x)^i}),\ i=2,1$ we have:
$$
X(2,z): y^2=x(1-x)(z-x),\ X(3,z): y^3=x(1-x)(z-x)^{2}.
$$
The genus of $X(2,z)$ (resp. $X(3,z)$) is $1$ (resp. $2$).
We make the transformation $x:=z-x$. Then
$$
X(6,z):  y^6=(z-x)(1-z+x)x^{5}.
$$
For $z=\frac{1}{2}$ we have the following new automorphism of $X(6,z)$
$$
\sigma(x,y)=(-x,\zeta_{12}y).
$$
\subsection{Decomposition of $T(6,z)$ into elliptic curves}
\begin{lemm} We have
\label{22aug06}
\begin{enumerate}
    \item $X(6,z)$ is hyperelliptic of genus 5, where
          \[ X(6,z) = \CC(x_1,y_1),\; y_1^2=x_1^{12}+(2-4z)x_1^6+1,\;
             y_1 =\frac{x^2-z^2+z}{x},\; x_1=\frac{y}{x}. \]
     \item \[ H^0(\CC(x_1,y_1),\Omega)=\langle x_1^i \frac{dx_1}{y_1} \mid i=0,\ldots, 4 \rangle \ni \frac{dx}{y}. \]
     \item If $z \neq 1/2$ then  $Aut(X(6,z))/\langle \epsilon \rangle=D_6=
           \langle \si,\tau \rangle,$ where
             $\epsilon$ denotes the hyperelliptic involution
            \[ \epsilon(x_1,y_1)=(x_1,-y_1),\quad
            \si (x_1,y_1)=(\zeta_6 x_1,y_1),\quad
               \tau (x_1,y_1)=( x_1^{-1},\frac{y_1}{x_1^6})\]
               resp.
           \[ \epsilon(x,y)=(-\frac{z-z^2}{x},-\frac{(z-z^2)y}{x^2}),  \quad
             \si (x,y)=(x,\zeta_6y), \quad
        \]
        \[
        \tau(x,y)=(\frac{(x-z)(1-z)}{1-z+x},\frac{x(x-z)(1-z)}{y(1-z+x)}). \]
       \item If $z=1/2$ then  $Aut(X(6,z))/\langle \epsilon \rangle=D_{12}=
           \langle \si,\tau \rangle,$ where
           \[ \si (x_1,y_1)=(\zeta_{12} x_1,y_1),\quad
               \tau (x_1,y_1)=( x_1^{-1},\frac{y_1}{x_1^6})\]
               resp.
           \[ \si (x,y)=(-x,\zeta_{12}y) \;
              \tau(x,y)=(\frac{(x-z)(1-z)}{1-z+x},\frac{x(x-z)(1-z)}{y(1-z+x)}). \]
   \item   We have \[ -6x_1^4 \frac{dx_1}{y_1}= \frac{dx}{y}=\eta_1.\]
\end{enumerate}
\end{lemm}
\begin{proof}
 \begin{enumerate}
    \item   Let $ K=X(6,z),\; y^6=(z-x)(1-z+x)x^5$,
           where $g(K)=5$ by the Hurwitz-Zeuthen genus formula.
    To each of the functions $ z-x,\; 1-z+x,\; x \in \CC(x)$ we associate
    the principal divisors
   $(z-x)= \frac{\p_z}{\n},\; (1-z+x)= \frac{\p_{z-1}}{\n},\;
      (x)= \frac{\p_0}{\n},\;$
    where the prime divisor $\n$ of $\CC(x)$
     corresponds to the place $p=\infty$
    and the  prime divisor $\p_0$ (resp. $\p_{z}$ and $\p_{z-1}$)
    corresponds to the prime polynomial $p=x$ (resp.
    $p=x-z$ and $p=x-z+1$) (see Chap. III,  \S 2 in \cite{Eich}).     
    The prime divisor $\p_0$ of $\CC(x)$  decomposes in $K$ into
    $\p_{0,1}^{e_1} \cdots \p_{0,r}^{e_r},$ where
    $\sum e_i=\deg(K/\CC(x))=6.$
    Comparing the order of the prime divisor $\p_{0,1}$ 
    in  the principal divisors
    $(y),(z-x),$ $(1-z+x),(x)$ of $K$ we get
    \[ \ord_{\p_{0,1}}(y^6)=6 \;\ord_{\p_{0,1}}(y)=\ord_{\p_{0,1}} ((z-x)(1-z+x)(x)^5)=    5 \;\ord_{\p_{0,1}}(x)= 5 e_1 .\]
    Hence $ \ord_{\p_{0,1}}(y) =5,\;  \ord_{\p_{0,1}}(x)=6$ and $r=1.$
    Thus the prime divisor $\p_0$ of $\CC(x)$ decomposes into
    $\p_{0,1}^6$ in $K$. (Hence we denote by  $\p_0$ also
    the prime divisor $\p_{0,1}$ of $K.$)
    
    Proceeding as above  we obtain that the principal divisors associated to 
           the functions $y,x,x-z,x-z+1$ in $K$ are
 \[ (y)=\frac{\p_0^5 \p_z \p_{z-1}}{\n^7},\; (x)=\frac{\p_0^6}{\n^6},\;
     (x-z)=\frac{\p_z^6}{\n^6},\; (x-z+1)=\frac{\p_{z-1}^6}{\n^6}. \]
Hence, the degree of the denominator of
\[ (\frac{y}{x})=\frac{\p_z \p_{z-1}}{\n \p_0} \]
is two. Thus $K$ has to be hyperelliptic, where
 $K=\CC(x_1,y_1),\;x_1=\frac{y}{x},\; y_1^2=f(x_1),\;
 f(x_1) \in \CC[x_1],\;\deg(f)=12$.
 Using that
\[ x_1^{6}= \frac{(z-x)(x-1+z)}{x} \]
and that there exists a $\si \in Aut(K/\CC)$ with $\si(x_1)=\zeta_6 x_1$
 we have to solve
 the following equation in order to determine $y_1:$
\[ x_1^{12}+a x_1^6+1 = y_1^2,\; a \in \CC. \]
Hence, we can assume
that $y_1= \frac{x^2+a_1x+a_2}{x},\; a_1,a_2 \in \CC.$
This gives after some computations
\[ a=2-4z,\; a_1=0, \; a_2=-z^2+z.\]

\item It is well known that
  \[ H^0(\CC(x_1,y_1),\Omega)=\langle x_1^i \frac{dx_1}{y_1} \mid i=0,\ldots, 4 \rangle. \]
  Since \[ (\frac{dx}{y})=\frac{{\mathfrak{d}}_{K/\CC(x)}}{\n_x^2} \frac{1}{(y)}=
               \frac{\p_0^5 \p_z^5 \p_{z-1}^5 \n^5}{\n^{12}}  \frac{\n^7}{\p_0^5 \p_z \p_{z-1}}=  \p_z^4 \p_{z-1}^4, \]
 where ${\mathfrak{d}}_{K/\CC(x)}$ denotes the pseudo different and
 ${\n_x}$ the divisor in $K$ of the denominator divisor of $x \in \CC(x)$,
(See Chap. III,  \S 2, \S3 in \cite{Eich}.)
   is an positive divisor the claim follows.
\end{enumerate}
Items $3,4$ and $5$ are proved by easy computations.
\end{proof}

We will show that  the vector space of holomorphic differentials
can be written as a direct sum $\oplus U_i,$ where $U_i=H^0(E_i,\Omega ), $
where $E_i$ is an elliptic subfield of $K.$
Thus we will study some subfields fixed by automorphisms of $Aut(X(6,z)):$

\begin{lemm}
\label{22aug}
 Let $\tau,\; \si \in Aut(K)$ be as in Lemma \ref{22aug06}.
 \[ H^0(X(6,z),\Omega )= \oplus_{i=1}^5 H^0(E_i,\Omega ), where \]
 \begin{enumerate}
   \item  $E_1$ is contained in the hyperelliptic function field
           $K^{\langle \si \rangle }$ of genus 2 and
        \[ E_1 : \tilde{y}^2=\tilde{x}^3-3\tilde{x}+2-4z,\;
              \tilde{x}=(x_1+x_1^{-1})^2-2,\; \tilde{y}=y_1x_1^{-3},\]
           where
     \[ j(E_1)=\frac{-2^4 3^3}{z(z-1)},\; \quad
 \frac{d \tilde{x}}{\tilde{y}}= 2(x_1^4-1)\frac{dx_1}{y_1}.\]

    \item  $E_2$ is contained in the hyperelliptic function field
           $K^{\langle \tau \si \rangle }$ of genus 2 and
           \[ E_2: \tilde{y}^2=\tilde{x}^3- 3\zeta_6^2\tilde{x}+2-4z,\;
           \tilde{x}=(x_1+\zeta_6 x_1^{-1})^2-2\zeta_6,\;
           \tilde{y}=y_1x_1^{-3},\] where
         \[    j(E_2)=\frac{-2^4 3^3}{z(z-1)},\; \quad \frac{d \tilde{x}}{\tilde{y}}=  2(x_1^4-\zeta_6^2)\frac{dx_1}{y_1}.\]

  \item  \[ E_3 =K^{\langle \si^2 \rangle}=X(2,z):  \tilde{y}^2= \tilde{x}^{4}+(2-4z)\tilde{x}^2+1,\;
                 \tilde{x}=x_1^3,\;  \tilde{y} = y_1, \] where
           \[ j(E_3)=2^8 \frac{(z^2-z+1)^3}{z^2(z-1)^2} \]
          and
 \[    \frac{d \tilde{x}}{\tilde{y}} =3 x_1^2  \frac{dx_1}{y_1}. \]
  \item   $E_4$ and $E_5$ are two isomorphic subfields
           of the hyperelliptic subfield $K^{\langle \si^3 \rangle}=X(3,z)$
                 of genus 2.
       \[ E_4 :  \tilde{y}^2=\tilde{x}^3-6\tilde{x}^2+9\tilde{x}-4z,\;
              \tilde{x}=x_1^4, \tilde{y}=y_1, \]
          where
           \[ j(E_4)=\frac{-2^4 3^3}{z(z-1)},\;\quad
          \frac{d \tilde{x}}{\tilde{y}}=2x_1 \frac{dx_1}{y_1} \]
          \[E_5:  \tilde{y}^2=\tilde{x}^3-6\tilde{x}^2+9\tilde{x}-4z,\;
                \tilde{x}=x_1^2+x_1^{-2},\; \tilde{y}=y_1x_1^{-6} \]
         where
             \[j(E_5)=\frac{-2^4 3^3}{z(z-1)},\; \quad
            \frac{d\tilde{x}}{\tilde{y}}=2x_1^3 \frac{dx_1}{y_1}.\]
 \end{enumerate}
\end{lemm}

\begin{proof}
 The claim follows from Lemma \ref{22aug06}.
\end{proof}

\begin{coro}
\label{elliptic} 
 Let $E_i,i=1,\ldots,5$ the elliptic subfields
  of $X(6,z)$ from Lemma \ref{22aug}.
 Then
 \[ H^0(X(6,z),\Omega)= \oplus_{i=1}^5 H^0(E_i,\Omega). \]
\end{coro}

\begin{proof}
 By Lemma \ref{22aug} we know that the elliptic differentials in $E_i,\;i=1,\ldots,5$ generate  $H^0(X(6,z),\Omega).$
\end{proof}

\begin{coro}
We have the isogeny
$$
T(6,z)\sim E^2,
$$
where
$$
E:\  y^2=x^3-3x+2-4z.
$$

\end{coro}
\begin{rem}\rm
 By a result of Ekedahl and Serre \cite{E-S} it follows already from the existence of the
 hyperelliptic subfield $X(3,z)$ of genus 2 that the Jacobian of $X(6,z)$
 is isogenous to a direct sum of elliptic curves.
Using the structure of the automorphism group it follows that
 $E_1$ is isogenous to $E_2$ by a result of Lange and Recillas \cite{L-R}.
\end{rem}
\subsection{ Proof of Theorem \ref{main1} }
\label{mainproof}
Recall the definition (\ref{22june2005}). 
Both the period sets  $\pedo(M_t,\nabla(\omega_0))$  
and $\pedo(E_t, \frac{dx}{y})$
satisfy the Gauss hypergeometric equation (\ref{gauss}) with
parameters $a,b,c$ given in (\ref{27.5.05}). This furnishes us with
the hypothesis of the Conjecture in \S \ref{gm}.
\begin{lemm}
\label{lemmagamma}
The space $\pedo(M_t,\nabla(\omega_0))$ is spanned over $\Q$ by
$\Gamma(\frac{1}{3})^3 \cdot \Q(\zeta_3)\cdot 
\pedo(E_t, \frac{dx}{y})$. 
\end{lemm}
\begin{proof}
Let $f_1=x_1^3-x_2^3-x_1-x_3$ and $f_2=t-(x_3^2+x_4^3+x_5^3)$ and fix a point
$b\in\C$ which is not a critical value neither for $f_1$ nor $f_2$.
Let also $\gamma:[0,1]\rightarrow \C, \ \gamma(\frac{1}{2})=b$ be a path in $\C$ which connects 
the critical point $t\in\C$ 
of $f_2$ to a critical point $c$ of $f_1$ such that $\delta_{1b}$ (resp. 
$\delta_{2b}$) vanishes along $\gamma$ (resp. $\gamma^{-1}$) when $s$ goes 
to $0$ (resp. $1$). One can show that the union
$$
\delta_t=\cup_{s\in [0,1]}\delta_{1\gamma(s)}\times \delta_{2\gamma(s)}
$$
defines a well-defined cycle in $H_4(M_t,\Z)$. In fact it is the join of
topological spaces $\delta_{1b}$ and $\delta_{2b}$ and
$$
I(t):=\int_{\delta_t} \omega_0=\int_0^1 I_1(s)I_2(s)ds,
$$
$$
I_1(s)=
\int_{\delta_{1\gamma(s)}}\frac{dx_1\wedge dx_2}{df_1},\ 
I_2(s)=\int_{\delta_{2\gamma(s)}}\frac{dx_3\wedge dx_4\wedge dx_5}{df_2}.
$$
(see \cite{arn} p. 53). A topological argument similar
to the one in \cite{arn} Theorem 2.1 or \cite{lam} \S 5 implies that
the cycles $\delta_t$ generate $H_4(M_t,\Q)$.
A simple integration shows that $I_2(s)$ is constant and  the space of such
integrals is equal to ${\Gamma(\frac{1}{3})^3} \Q(\zeta_3)$ (see \cite{dmos} 
Lemma 7.12). 
This means that up to  $ {\Gamma(\frac{1}{3})^3}\Q(\zeta_3) $ the integral 
$I(t)$ reduces to $\int_\Delta dx_1\wedge dx_2$,  
where $\Delta:=\cup_{s\in [0,1]}\delta_{1\gamma(s)}\in H_2(\C^2, f_1^{-1}(t),\Z)$. By Stokes theorem this  integral is 
$$
\int_{\delta_{1t}}x_1dx_2.
$$
Since $\frac{\partial I}{\partial t}=\int_{\delta_t}\nabla(\omega_0)$, we 
conclude that
$$
\int_{\delta_t}\nabla(\omega_0)\in {\Gamma(\frac{1}{3})^3} \Q(\zeta_3)\pedo 
(\{f_1=t\}, \frac{dx_1\wedge dx_2}{df_1}).
$$  
Now the  elliptic curve $\{f_1=t\}$ is isomorphic to 
the one in (\ref{26.5}) by the mapping 
$$
(x_1,x_2)\rightarrow (\frac{3t}{x_1+x_2}+1, \ 9t( \frac{ x_2}{ x_1+x_2}-
\frac{1}{2}))
$$
and under this mapping  $\frac{dx_1\wedge dx_2}{df_1}$ goes to $-\frac{dx}{2y}$
\end{proof}

Now let us prove Theorem \ref{main1}.   The proof of
the first and second part is as follows: 
The map $\delta\in H_4(M_t,\Q(\zeta_3))
\rightarrow \int_{\delta}\nabla(\omega_0)\in 
\pedo(M_t,\nabla(\omega_0))$ is surjective. Now  Lemma \ref{lemmagamma}, Corollary
\ref{bewerbung} and Corollary \ref{elliptic} 
imply that $\pedo(M_t,\nabla(\omega_0))$ is of dimension one if and only if
 $E_t$ is elliptic.

By Proposition \ref{deligne} if
$\delta\in H_4(M_t,\Q)$ is a Hodge cycle then
\begin{equation}
\label{footbal}
\frac{\partial}{\partial t}(\int_{\delta}\nabla(\omega_0))=\int_{\delta}\nabla^2(\omega_0)\in 
\pi^2\bar\Q.
\end{equation}
Let $Y_1(t),Y_2(t)$ be a basis of 
$\pedo(E_t,\frac{dx}{y})$ and $\frac{Y_1(t_0)}{Y_2(t_0)}=
r\in\Q(\zeta_3)$ at some point $t_0$. Then according to
Lemma \ref{lemmagamma} we have $\Gamma(\frac{1}{3})^3Y_1(t)=
\int_{\delta_{1,t}}\nabla(\omega_0)$ and $\Gamma(\frac{1}{3})^3rY_2(t)=
\int_{\delta_{2,t}}\nabla(\omega_0)$ for some $\delta_{1,t},\delta_{2,t} 
\in H_4(M_t,\Q)$ 
and all $t$ in a neighborhood of $t_0$.
This implies that $\delta_{1,t_0}-\delta_{2,t_0}$ is a Hodge cycle. 
Now we use  (\ref{footbal}) and we conclude that 
$$
\Gamma(\frac{1}{3})^3(Y_1'(t_0)-rY_2'(t_0))\sim \pi^2.
$$
By Legendre theorem on the determinant of the period matrix 
of elliptic curves (or the relation (\ref{13july05})),
we have 
$$
Y_1'(t_0)-aY_2'(t_0)=\frac{(Y_1'Y_2-Y_1Y_2')(t_0)}{Y_2(t_0)}\sim \frac{\pi}{Y_2(t_0)}
$$
and so $Y_1(t_0),Y_2(t_0)\sim \frac{1}{\pi}\Gamma(\frac{1}{3})^3$. 
Using the formula
(\ref{akh}) the third part of the above theorem is proved. 
\section{Other examples}
In this section we give a list of cubic hypersurfaces such that
the corresponding integral $\int_{\delta_t} \nabla(\frac{dx}{f})$ 
for the definition of Hodge cycles satisfy
a Picard-Fuchs equation of order two. For calculations 
we have used the library {\tt brho.lib} written in {\sc Singular} 
(see \cite{hos005, mo}).   Let $\omega_0=\frac{dx}{df}$.
\begin{enumerate}
\item  
For $$f=x_1^3+\cdots+x_5^3-x_1x_2
$$ 
$\int_{\delta_t} \omega_0$  satisfies 
$$
(27t^2+t)y''+6y=0
$$
and so 
$\int_{\delta_t}\nabla\omega_0$ satisfies
$$
(27t^2+t)y''+54ty'-6y=0.
$$
\item
For
$$
f=x_1^3+x_2^3+x_3^3+x_4^3+x_5^3-x_1^2-x_2^2 
$$
$\int_{\delta_t} \omega_0 $ satisfies
$$
(405t+60)y'+
(2187t^2+648t+32)y''+
(729t^3+324t^2+32t)y'''=0
$$
and $y'= \int_{\delta_t} \nabla \omega_0$. Note that
the coefficient of $y'''$ is $t(27t+4)(27t+8)$.
\item
For
$$
f=x_1^3+x_2^3+x_3^3+x_4^3+x_5^3-x_1^2-x_2x_1
$$
$\int_{\delta_t} \omega_0 $ satisfies
$$
0y+
(4374t^2-1296t-108)y'+
(39366t^3+3645t^2-1188t+2)y''+
$$
$$
(19683t^4+6561t^3-621t^2+2t)y'''=0.
$$
The coefficient of $y''$ is $t(729t^2+297t-1)(27t-2)$ and so there is an 
apparent singularity at
$t=\frac{2}{27}$.  Note that the critical values of $f$ are the roots 
of $(729t^2+297t-1)t$.
\item
For
$$
f=x_1^3+x_2^3+x_3^3+x_4^3+x_5^3-x_1^2-x_1
$$
$\int_{\delta_t} \omega_0$ satisfies
$$
(27t+11)y'+
(81t^2+66t-15)y''=0.
$$
The coefficient of $y''$ is $3(27t-5)(t+1)$. 
\item
For
$$
f=x_1^3+x_2^3+x_3^3+x_4^3+x_5^3-x_1-x_1x_2
$$
the set of critical values is the roots of 
$S(t)=(19683t^4+2187t^3-5751t^2+541t+433)$ and
$\int_{\delta_t} \omega_0$ satisfies
$$
(3188646t^4+472392t^3+6403536t^2+236844t-314922)y'+
$$
$$
(28697814t^5+4782969t^4+25824096t^3+1427382t^2-4369464t+152443)y''+
$$
$$
(14348907t^6+2657205t^5+2322594t^4+794610t^3-1524204t^2+199207t+140725)
y'''=0.
$$
The coefficient of $y'''$ is
$$
(19683t^4+2187t^3-5751t^2+541t+433)(729t^2+54t+325)
$$
and so there are two apparent singularities.

\end{enumerate}

\def\cprime{$'$} \def\cprime{$'$}

\bibliographystyle{plain}

\end{document}